\documentclass[11pt]{article}
\usepackage{enumerate}
\usepackage{amssymb,a4wide,latexsym,makeidx,epsfig,fleqn}
\usepackage{amsthm}
\usepackage{amsmath}
\usepackage{enumerate}
\usepackage{graphicx}
\usepackage{epstopdf}
\newtheorem{theorem}{Theorem}[section]

\newtheorem{lemma}[theorem]{Lemma}

\begin{document}
\textwidth 150mm \textheight 225mm
\title{Distance signless Laplacian spectral radius and tough graphs involving minimun degree\footnote{This work is supported by the National Natural Science Foundations of China (No. 12371348, 12201258), the Postgraduate Research \& Practice Innovation Program of Jiangsu Normal University (No. 2024XKT1702).}}
\author{{Xiangge Liu, Yong Lu\footnote{Corresponding author}, Caili Jia, Qiannan Zhou, Yue Cui}\\
{\small  School of Mathematics and Statistics, Jiangsu Normal University,}\\ {\small  Xuzhou, Jiangsu 221116,
People's Republic
of China.}\\
{\small E-mail:  luyong@jsnu.edu.cn}}

\date{}
\maketitle
\begin{center}
\begin{minipage}{120mm}
\vskip 0.3cm
\begin{center}
{\small {\bf Abstract}}
\end{center}
{\small Let $G=(V(G),E(G))$ be a simple graph, where $V(G)$ and $E(G)$ are the vertex set and the edge set of $G$, respectively.
The number of components of $G$ is denoted by $c(G)$.
Let $t$ be a positive real number, and a connected graph $G$ is $t$-tough if $t c(G-S)\leq|S|$ for every vertex cut $S$ of $V(G)$.
The toughness of graph $G$, denoted by $\tau(G)$, is the largest value of $t$ for which $G$ is $t$-tough.
Recently, Fan, Lin and Lu [European J. Combin. 110(2023), 103701] presented sufficient conditions based on the spectral radius for  graphs to be  1-tough with minimum degree $\delta(G)$ and  graphs to be $t$-tough with $t\geq 1$ being an integer, respectively.
In this paper, we establish sufficient conditions in terms of the distance signless Laplacian spectral radius for graphs to be 1-tough with minimum degree $\delta(G)$ and  graphs to be $t$-tough, where $\frac{1}{t}$ is a positive integer.
Moreover, we consider the relationship  between the distance signless Laplacian spectral radius and $t$-tough graphs in terms of the order $n$.

 }
\end{minipage}
\end{center}

\section{Introduction }
\hspace{1.3em}
Let $G=(V(G),E(G))$ be a simple graph, where $V(G)$ is the vertex set and $E(G)$ is the edge set.
The order and size of $G$ are denoted by $|V(G)|=n$ and $|E(G)|=m$, respectively.
For a vertex subset $S$ of $G$, we denote by $G-S$ and $G[S]$ the subgraph of $G$ obtained from $G$ by deleting the vertices in $S$ together with their incident edges and the subgraph of $G$ induced by $S$, respectively.
The number of components of $G$ is denoted by $c(G)$.
Let $t$ be a positive real number and a connected graph $G$ is \emph{$t$-tough} if $tc(G-S)\leq|S|$ for every vertex cut $S$ of $V(G)$.
The\emph{ toughness} of graph $G$, denoted by $\tau(G)$, is the largest value of $t$ for which $G$ is $t$-tough (taking $\tau(K_{n})=\infty$, where $K_{n}$ is a \emph{complete graph} of order $n$).
Thus, $\tau(G)=\min\{\frac{|S|}{c(G-S)}:S\subseteq V(G),c(G-S)\geq2\}$.
The concept of toughness was initially proposed by Chv\'{a}tal \cite{C} in 1973, which serves as a simple way to measure how tightly various pieces of a graph hold together.
Let $d(v)$ be the \emph{degree} of vertex $v\in V(G)$, $\delta(G)$ be the \emph{minimum degree} ($\delta$ for short) of $G$.
For two vertex-disjoint graphs $G_{1}$ and $G_{2}$, we use $G_{1}+G_{2}$ to denote the \emph{disjoint union} of $G_{1}$ and $G_{2}$.
The \emph{join} $G_{1}\vee G_{2}$ is the graph obtained from $G_{1}+G_{2}$ by adding all possible edges between $V(G_{1})$ and $V(G_{2})$.
We use $J$ to denote the all-one matrix, $I$ to denote the identity square matrix.

For a simple graph $G$ of order $n$,
the \emph{adjacency matrix} of $G$ is denoted by $A(G)$, and $\lambda_{1}(G)\geq\lambda_{2}(G)\geq\cdots\geq\lambda_{n}(G)$ are the eigenvalues of $A(G)$.
In particular, the eigenvalue $\lambda_{1}(G)$ is called the \emph{spectral radius} of $G$ (short for $\rho(G)$).
The \emph{distance} between $v_{i}$ and $v_{j}$ denoted by $d_{ij}(G)$, is the length of a shortest path from $v_{i}$ to $v_{j}$, and the \emph{distance matrix} of $G$ is the symmetric matrix $D(G)=(d_{ij}(G))_{n\times n}$  with its rows and columns indexed by $V(G)$.
The \emph{transmission} $Tr(v)$ of a vertex $v \in V(G)$ is the sum of distances from $v$ to all vertices in $G$.
We say that a graph is \emph{$k$-transmission-regular} (or transmission-regular) if its distance matrix has constant row sum equal to $k$.
Aouchiche and Hansen \cite{AHP} gave the definition of the \emph{distance signless Laplacian matrix} of a graph $G$, which is given by $Q_{D}(G) =D(G)+Tr(G)$, where $Tr(G)$ is the diagonal matrix whose diagonal entries are the vertex transmissions in $G$.
The largest eigenvalue $\eta_{1}(G)$ is called the \emph{distance signless Laplacian spectral radius} of graph $G$.

The connections between the structural properties and  the distance signless Laplacian spectral radius have been well studied.
Zhou and Wang \cite{ZW} mainly considered the relationships between distance signless Laplacian spectral radius and the Hamiltonian properties of graphs according to the size $m$ and the order $n$.
Zhou et al. \cite{ZZL} investigated the relations between the spanning $k$-tree and the distance signless Laplacian spectral radius in a connected graph and proved an upper bound for $\eta_{1}(G)$ in a connected graph $G$ to guarantee the existence of a spanning $k$-tree.
Mahato and Kannan \cite{MK} showed that the generalized tree shift transformation increases the distance spectral radius and the distance signless Laplacian spectral radius of complements of trees.

Many researchers have paid attention to the problems of establishing relations between  the toughness and  spectra of graphs.
In 2023, Fan et al. \cite{FLL} presented sufficient conditions based on the spectral radius for  graphs to be  1-tough with minimum degree $\delta(G)$ and graphs to be $t$-tough with $t\geq 1$ being an integer, respectively.
Chen et al. \cite{CFL} considered a new parameter $\tau(G)$ which is a slight variation of toughness and it was defined as $\tau'(G)=\min\{\frac{|S|}{c(G-S)-1}:S\subseteq V(G),c(G-S)\geq2\}$ by Enomoto \cite{EJK}.
They investigated an analogous problem concerning balanced bipartite graphs and provided spectral
radius and edge conditions for 2-factors in 1-tough balanced bipartite graphs by incorporating the toughness and spectral conditions.
Recently, Lou et al. \cite{LLS} presented a sufficient condition based on the distance spectral radius to guarantee that a graph is 1-tough with minimum degree $\delta(G)$.
The \emph{signless Laplacian matrix} of $G$ is defined as $Q(G)=D'(G)+A(G)$, where $D'(G)$ is the diagonal matrix of vertex degrees of $G$, and the largest eigenvalue of $Q(G)$ is called the \emph{signless Laplacian spectral radius} of $G$.
Moreover, Jia and Lou \cite{JL}  considered the signless Laplacian spectral radius versions of the  above problem (a graph is 1-tough with minimum degree $\delta(G)$).  For more extensive work on toughness related to  eigenvalues of graphs, one can see $(\cite{CL},\cite{CW},\cite{HDZ})$.


In this paper, we consider a sufficient condition in terms of the distance signless Laplacian spectral radius $\eta_{1}(G)$ to ensure that a graph $G$ is 1-tough with minimum degree $\delta$.

\noindent\begin{theorem}\label{th:1.1.}
Let $G$ be a connected graph of order $n\geq\max\{11\delta,\frac{1}{2}\delta^{2}+2\delta\}$ with minimum degree $\delta\geq2$.
If $$\eta_{1}(G)\leq\eta_{1}(K_{\delta}\vee(K_{n-2\delta}+\delta K_{1})),$$ then $G$ is 1-tough unless $G\cong K_{\delta}\vee(K_{n-2\delta}+\delta K_{1})$.
\end{theorem}

For clarity, the exceptional graph in Theorem  \ref{th:1.1.} as shown in Figure 1.

\begin{figure}[htbp]
  \centering
 \includegraphics[scale=0.9]{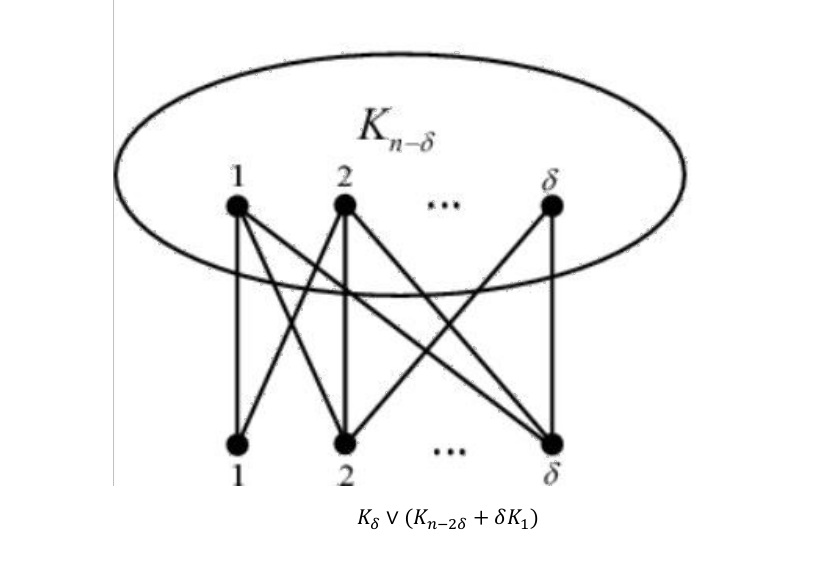}
 \caption{$K_{\delta}\vee(K_{n-2\delta}+\delta K_{1})$ .}
\end{figure}

Chen et al. \cite{CLW} further showed sufficient eigenvalue conditions
for a regular graph and  a $k$-regular bipartite graph to be $\frac{1}{t}$-tough for a positive integer $t$.
Additionally, Lou et al. \cite{LLS}  presented a sufficient condition based on the distance spectral radius to guarantee that a graph is $t$-tough, where $t$ or $\frac{1}{t}$ is a positive integer.
 Jia and Lou \cite{JL} also considered the above problem  based on the signless Laplacian spectral radius.

Motivated by their results, we are devoted to giving a sufficient condition in terms of the distance signless Laplacian spectral radius $\eta_{1}(G)$ for graphs to be $t$-tough, where $\frac{1}{t}$ is a positive integer.

\noindent\begin{theorem}\label{th:1.2.}
Let $G$ be a connected graph of order $n\geq \frac{2}{t^{2}}+\frac{3}{t}+3t+4$ and $\frac{1}{t}$ be a positive integer.
If $\eta_{1}(G)\leq\eta_{1}(K_{1}\vee(K_{n-\frac{1}{t}-1}+\frac{1}{t}K_{1}))$, then $G$ is $t$-tough unless $G\cong K_{1}\vee(K_{n-\frac{1}{t}-1}+\frac{1}{t}K_{1})$.
\end{theorem}

For clarity, the exceptional graph in Theorem \ref{th:1.2.} as shown in Figure 2.

\begin{figure}[htbp]
  \centering
 \includegraphics[scale=0.9]{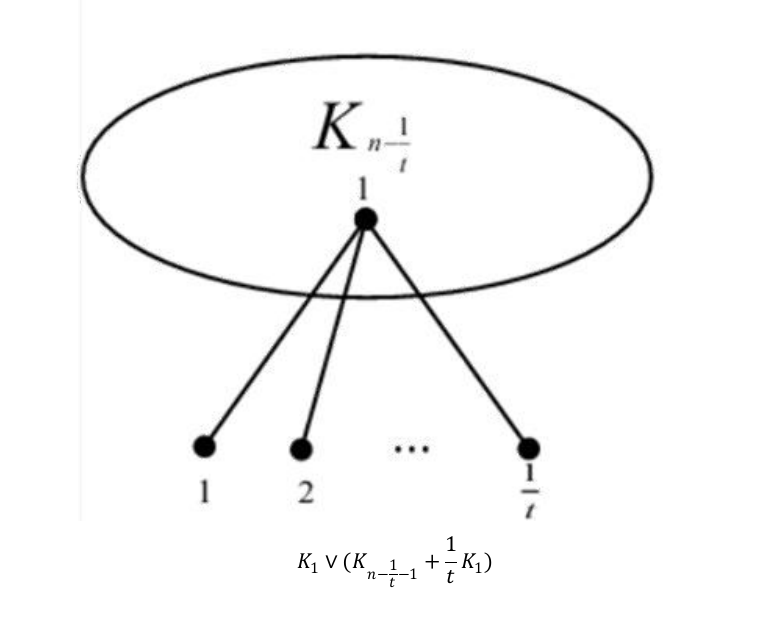}
 \caption{ $K_{1}\vee(K_{n-\frac{1}{t}-1}+\frac{1}{t}K_{1})$.}
\end{figure}

Bondy \cite{bondy} conjectured that almost any nontrivial condition on a graph which implies that the graph is Hamiltonian also implies that the graph is pancyclic.
Benediktovich \cite{B} confirmed the Bondy conjecture for $t$-tough graphs in the case when $t
\in\{1,2,3\}$ in terms of the size, the spectral radius and the signless Laplacian spectral radius.
Chen et al. \cite{CLX} considered the new parameter $\tau'(G)$ which is mentioned above and established two lower bounds on the size to guarantee that a graph $G$ is $t$-tough.

In this paper, we also present  the relationships  between the distance signless Laplacian spectral radius and graphs to be $t$-tough in terms of the order $n$ according to the new parameter $\tau'(G)=\min\{\frac{|S|}{c(G-S)-1}\}$.

\noindent\begin{theorem}\label{th:1.3.}
Let $G$ be a connected graph of order $n$ and minimum degree $\delta$.
Then the following statements hold.
\begin{enumerate}[(a)]
\item
Let $t\geq2$ and $n\geq2t^{2}+2t$ be two integers.
If $\eta_{1}(G)<\frac{2n^{2}+2n-4t}{n}$, then $G$ is $t$-tough.
\item
Let $n\geq\max\{(\frac{2t^{2}+3t+1}{t})\delta+\frac{1}{2t}+5,\frac{1}{6}(\frac{\delta(\delta+4)+1}{t}+9\delta+6)\}$, where $\frac{1}{t}\geq1$ and $\delta\geq t+1$ are integers.
If $\eta_{1}(G)<2n+\frac{4\delta}{t}+2-\frac{2(\delta+t)(2\delta t+\delta+2t)}{nt^{2}}$, then $G$ is $t$-tough.
\end{enumerate}
\end{theorem}

\section{Proof of Theorem \ref{th:1.1.} }
\hspace{1.3em}
In this section, we first present a fundamental result to compare the distance signless Laplacian spectral radius of a graph and its spanning subgraph.

\noindent\begin{lemma}\label{le:2.1.}\cite{MH}
Let $e$ be an edge of a graph $G$ such that $G-e$ is still connected. Then $\eta_{1}(G-e)>\eta_{1}(G)$.
\end{lemma}

Next, we explain the concepts of \emph{equitable matrices} and \emph{equitable partitions} in \cite{BH}.
Let $M$ be a real $n\times n$ matrix  and let $X=\{1,2,\ldots,n\}$.
Given a partition  $\Pi=\{X_{1}, X_{2}, \ldots, X_{m}\}$ with $X=X_{1}\cup X_{2}\cup\cdots \cup X_{m}$, the matrix can be described in the following block form
\begin{align*}
M=\left(
\begin{array}{ccccccccc}
M_{11} & M_{12} & \cdots & M_{1m} \\
M_{21} & M_{22} & \cdots & M_{2m} \\
\vdots & \vdots & \ddots &\vdots \\
M_{m1} & M_{m2} & \cdots & M_{mm} \\
\end{array}
\right).
\end{align*}
The \emph{quotient matrix} $R(M)$ of the matrix $M$ (with respect to the given partition) is the $m\times m$ matrix whose entries are the average row sums of the blocks $M_{i,j}$ of $M$.
The above partition is called \emph{equitable} if each block $M_{i,j}$ of $M$ has constant row sum.

\noindent\begin{lemma}\label{le:2.2.}\cite{YYSX}
Let $R(M)$ be an equitable quotient matrix of $M$ as defined above, and $M$ be a nonnegative matrix.
Then $\rho(R(M))=\rho(M)$, where $\rho(R(M))$ and $\rho(M)$ denote the largest eigenvalues of the matrices $R(M)$ and $M$.
\end{lemma}

The \emph{Wiener index} $W(G)$ of a connected graph $G$ of order $n$ is defined by the sum of all distances in $G$, that is, $W(G)=\sum\limits_{i<j}d_{ij}(G)$.
Xing et al. \cite{XZL} presented a lower bound on the distance signless Laplacian spectral radius of a graph.

\noindent\begin{lemma}\label{le:2.3.}\cite{XZL}
Let $G$ be a connected graph of order $n$.
Then $$\eta_{1}(G)\geq\frac{4W(G)}{n}$$
with equality if and only if $G$ is transmission regular.
\end{lemma}

\noindent\begin{lemma}\label{le:2.4.}\cite{LL}
Let $n,q,s$ and $n_{i}(i=1,2,\ldots,c)$ be positive integers with $n_{1}\geq n_{2}\geq\cdots\geq n_{c}\geq1$ and $n=s+\sum\limits_{i=1}^{c}n_{i}$.
Then
$$\eta_{1}(K_{s}\vee (K_{n_{1}}+K_{n_{2}}+\cdots+K_{n_{c}}))\geq\eta_{1}(K_{s}\vee (K_{n-s-(c-1)}+(c-1)K_{1}))$$ with equality if and only if $K_{s}\vee (K_{n_{1}}+K_{n_{2}}+\cdots+K_{n_{c}})\cong K_{s}\vee (K_{n-s-(c-1)}+(c-1)K_{1})$.
\end{lemma}

According to Lemma \ref{le:2.4.}, we give the following critical lemma which generalizes the above lemma.

\noindent\begin{lemma}\label{le:2.5.}
Let $n,q,s$ and $n_{i}(i=1,2,\ldots,c)$ be positive integers with $n_{1}-n_{2}\geq 2p$, $n_{1}\geq n_{2}\geq\cdots\geq n_{c}\geq p\geq1$ and $n=s+\sum\limits_{i=1}^{c}n_{i}$.
Then
$$\eta_{1}(K_{s}\vee (K_{n_{1}}+K_{n_{2}}+\cdots+K_{n_{c}}))\geq\eta_{1}(K_{s}\vee (K_{n-s-(c-1)p}+(c-1)K_{p}))$$
with equality if and only if $K_{s}\vee (K_{n_{1}}+K_{n_{2}}+\cdots+K_{n_{c}})\cong K_{s}\vee (K_{n-s-(c-1)p}+(c-1)K_{p})$.
\end{lemma}
\noindent\textbf{Proof.}

Let $G_{1}=K_{s}\vee (K_{n-s-(c-1)p}+(c-1)K_{p})$.
The distance signless Laplacian matrix of $G_{1}$ is $Q_{D}(G_{1})=D(G_{1})+Tr(G_{1})$, where
\begin{align*}
\scalebox{0.8}{$
D(G_{1})=
\bordermatrix{%
& s & n-s-(c-1)p & p & \cdots & p\cr
s & J-I & J & J & \cdots & J \cr
n-s-(c-1)p & J & J-I & 2J & \cdots & 2J \cr
p & J & 2J & J-I & \cdots & 2J \cr
\vdots & \vdots & \vdots & \vdots & \ddots &\vdots \cr
p & J & 2J & 2J & \cdots & J-I \cr
},$}
\end{align*}
\begin{align*}
\scalebox{0.8}{$
Tr(G_{1})=
\bordermatrix{%
& s & n-s-(c-1)p & p & \cdots & p\cr
s & (n-1)I & 0 & 0 & \cdots & 0 \cr
n-s-(c-1)p & 0 & (n+(c-1)p-1)I & 0 & \cdots & 0 \cr
p & 0 & 0 & (2n-s-p-1)I & \cdots & 0 \cr
\vdots & \vdots & \vdots & \vdots & \ddots &\vdots \cr
p & 0 & 0 & 0 & \cdots & (2n-s-p-1)I \cr
}.$}
\end{align*}
The identity square matrix $I$ in $(n-1)I$ is the $s\times s$ matrix, and in $(n+(c-1)p-1)I$ is the $(n-s-(c-1)p)\times (n-s-(c-1)p)$ matrix.

Hence,
the equitable quotient matrix of the distance signless Laplacian matrix $Q_{D}(G_{1})$ is equal to
\begin{align*}
B=
\left(\begin{matrix}
n+s-2 & n-s-(c-1)p & (c-1)p  \cr
s & 2n-s-2 & 2(c-1)p  \cr
 s & 2\big(n-s-(c-1)p\big) & 2n-s+2(c-2)p-2  \cr
\end{matrix}\right).
\end{align*}

By the Perron-Frobenius theorem, $\eta_{1}(G)$ is always positive (unless $G$ is a trivial graph which has one vertex).
Furthermore, when $G$ is connected, there exists a unique positive unit eigenvector $\textbf{x}=(x_{1}, x_{2}, \ldots, x_{n})^{T}$ corresponding to $\eta_{1}(G)$, which is called a \emph{Perron vector} of $Q_{D}(G)$.

Let $\textbf{x}$ be the Perron vector of $Q_{D}(G_{1})$ corresponding to $\eta_{1}(G)$, and we take $x_{u}=x_{1}$ for all $u\in V(K_{s})$, $x_{v}=x_{2}$ for all $x_{v}\in V(K_{n-s-(c-1)p})$, and $x_{w}=x_{3}$ for all $x_{w}\in V((c-1)K_{p})$.
Due to $Q_{D}(G_{1})\textbf{x}=\eta_{1}(G_{1})\textbf{x}$,
$$\left\{
\begin{aligned}
&\eta_{1}(G_{1})x_{2}=sx_{1}+(2n-s-2)x_{2}+2(c-1)px_{3}\\
&\eta_{1}(G_{1})x_{3}=sx_{1}+2\big(n-s-(c-1)p\big)x_{2}+\big(2n-s+2(c-2)p-2\big)px_{3}.
\end{aligned}
\right.$$
Thus,
$x_{3}=\frac{\eta_{1}(G_{1})-s-2(c-1)p+2}{\eta_{1}(G_{1})-2n+s+2p+2}x_{2}$,
where $\eta_{1}(G_{1})>0$ and $p\geq1$, then $$x_{2}+x_{3}=\frac{2\big(\eta_{1}(G_{1})+2-n-(c-2)p\big)}{\eta_{1}(G_{1})-2n+s+2p+2}x_{2}.$$

Let $G_{2}=K_{s}\vee(K_{n_{1}}+K_{n_{2}}+\cdots+K_{n_{c}})$. The distance signless Laplacian matrix of $G_{2}$ is $Q_{D}(G_{2})=D(G_{2})+Tr(G_{2})$, where
\begin{align*}
D(G_{2})=
\scalebox{0.9}{$
\bordermatrix{%
& s & n_{1} & n_{2} & \cdots & n_{c}\cr
s & J-I & J & J & \cdots & J \cr
n_{1} & J & J-I & 2J & \cdots & 2J \cr
n_{2} & J & 2J & J-I & \cdots & 2J \cr
\vdots & \vdots & \vdots & \vdots & \ddots &\vdots \cr
n_{c} & J & 2J & 2J & \cdots & J-I \cr
},$}
\end{align*}
\begin{align*}
\scalebox{0.9}{$
Tr(G_{2})=
\bordermatrix{%
& s & n_{1} & n_{2} & \cdots & n_{c}\cr
s & (n-1)I & 0 & 0 & \cdots & 0 \cr
n_{1} & 0 & (2n-n_{1}-s-1)I & 0 & \cdots & 0 \cr
n_{2} & 0 & 0 & (2n-n_{2}-s-1)I & \cdots & 0 \cr
\vdots & \vdots & \vdots & \vdots & \ddots &\vdots \cr
n_{c} & 0 & 0 & 0 & \cdots & (2n-n_{c}-s-1)I \cr
}.$}
\end{align*}

Let $\textbf{\~{x}}=(\underbrace{x_{1}, \ldots, x_{1}}_{s}, \underbrace{x_{3}, \ldots, x_{3}}_{(c-1)p}, \underbrace{x_{2}, \ldots, x_{2}}_{n-s-(c-1)p})^{T}$.
The matrix $Q_{D}(G_{2})-Q_{D}(G_{1})$ is as follows:

\begin{align*}
\scalebox{0.8}{$
\bordermatrix{%
& s & (c-1)p & n_{1} & n_{2}-p & \cdots & n_{c}-p\cr
s      & 0 & 0 & 0 & 0 & \cdots & 0 \cr
(c-1)p & 0 & D_{1} & 0 & -E_{1} & \cdots & -E_{c-1} \cr
n_{1}  & 0 & 0 & (n-n_{1}-(c-1)p-s)I & J & \cdots & J \cr
n_{2}-p& 0 & -E_{1}^{T} & J & (n-n_{2}-(c-1)p-s)I &  \cdots & J\cr
\vdots & \vdots & \vdots & \vdots & \vdots &  & \ddots &\vdots \cr
n_{c}-p & 0 & -E_{c-1}^{T} & J & J  & \cdots & (n-n_{c}-(c-1)p-s)I\cr
},$}
\end{align*}
where $E_{i}$ denotes the matrix whose each entry from $[(i-1)p+1]$-th row to $ip$-th row is 1 and whose other entries are 0.
Moreover, $D_{1}=diag\{d_{2},d_{3},\ldots,d_{c}\}$, where $d_{i}=(p-n_{i})I_{p \times p}$.


Thus,
\begin{align}
&\eta_{1}(G_{2})-\eta_{1}(G_{1})\geq \textbf{\~{x}}^{T}(Q_{D}(G_{2})-Q_{D}(G_{1}))\textbf{\~{x}}\nonumber
\\=&n_{1}\sum_{i=2}^{c}(n_{i}-p)x_{2}^{2}+\sum_{i=2}^{c}(n_{i}-p)\big[\big(n-s-n_{i}-(c-1)p\big)x_{2}^{2}-2px_{2}x_{3}\big]\nonumber
\\&+n_{1}\big[n-s-n_{1}-(c-1)p\big]x_{2}^{2}+\sum_{i=2}^{c}p(p-n_{i})x_{3}^{2}\nonumber
\\&+\sum_{i=2}^{c}(n_{i}-p)\big[n-s-n_{i}-(c-1)p\big]x_{2}^{2}\nonumber
\\=&\sum_{i=2}^{c}(n_{i}-p)\Big[\big(2n_{1}+2n-2s-2n_{i}-2(c-2)p\big)x_{2}^{2}-px_{2}^{2}-2px_{2}x_{3}-px_{3}^{2}\Big]\nonumber
\\=&\sum_{i=2}^{c}(n_{i}-p)\Big[\big(2n_{1}+2n-2s-2n_{i}-2(c-2)p\big)x_{2}^{2}-p(x_{2}+x_{3})^{2}\Big]\nonumber
\\=&2x_{2}^{2}\sum_{i=2}^{c}(n_{i}-p)\Big[\big(n_{1}+n-s-n_{i}-(c-2)p\big)-2p(\frac{\eta_{1}(G_{1})+2-n-(c-2)p}{\eta_{1}(G_{1})-2n+s+2p+2})^{2}\Big].\nonumber
\end{align}

Note that $n_{1}\geq n_{2}\geq\cdots\geq n_{c}\geq p$ and $n_{1}-n_{2}\geq2p$.
In order to prove $\eta_{1}(G_{2})-\eta_{1}(G_{1})>0$, we need to prove $$\big(n_{1}+n-s-n_{i}-(c-2)p\big)-2p\big(\frac{\eta_{1}(G_{1})+2-n-(c-2)p}{\eta_{1}(G_{1})-2n+s+2p+2}\big)^{2}>0,$$ which implies $\eta_{1}(G_{1})>2n-s-2p-2+r$, where $r=\frac{n-s-cp}{\sqrt{\frac{n-s-(c-4)p}{2p}}-1}$.

Consider the determinant of matrix $(2n-s-2p-2+r)I-B$, if $det((2n-s-2p-2+r)I-B)<0$, we can get $\eta_{1}(G_{1})=\rho(B)>2n-s-2p-2+r$.
Thus,

\begin{align*}
det&\big((2n-s-2p-2+r)I-B\big)=
\begin{vmatrix}
n-2s-2p+r & -n+s+(c-1)p & -(c-1)p \\
-s & r-2p & -2(c-1)p \\
-s & -2n+2s+2(c-1)p & r-2(c-1)p\\
\end{vmatrix}\nonumber
\\=&r
\begin{vmatrix}
n-2s-2p+r & -n+s+(c-1)p  \\
-s & r-2p  \\
\end{vmatrix}
+2(c-1)p
\begin{vmatrix}
n-\frac{3}{2}s-2p+r & -n+s+cp-\frac{1}{2}r \\
0 & -2n+2s+2cp-r\\
\end{vmatrix}
\\=&r\Big[(r-2p)(n-2s-2p+r)+s(s+(c-1)p-n)\Big]\nonumber
\\&-2(c-1)p(n-\frac{3}{2}s-2p+r)(2n-2s-2cp+r)\nonumber
\\=&r\Big[(r-2p)(n-2s-2p)+r^{2}-2rp-s(n-s-(c-1)p)\Big]\nonumber
\\&-2(c-1)p(n-\frac{3}{2}s-2p+r)(2n-2s-2cp+r)\nonumber
\\=&\Big[(n-2s-4p)r^{2}-2(c-1)pr^{2}-4(c-1)p(n-\frac{3}{2}s-2p)(n-s-cp)\Big]\nonumber
\\&-2p(n-2s-2p)r-prs-2(c-1)p(n-\frac{3}{2}s-2p)r\nonumber
\\&+\Big[r^{2}-(n-s-cp)s-4(c-1)p(n-s-cp)\Big]r.\nonumber
\end{align*}

Hence, we need to check whether inequalities \begin{align}
(n-2s-4p)r^{2}-2(c-1)pr^{2}-4(c-1)p(n-\frac{3}{2}s-2p)(n-s-cp)<0
\end{align}
and
\begin{align}
r^{2}-(n-s-cp)s-4(c-1)p(n-s-cp)<0
\end{align}
hold.

If $n-2s-2(c+1)p\leq0$, then the inequality (1) holds.
If $n-2s-2(c+1)p>0$, then the inequality (1) implies
\begin{align}
&\big(n-2s-2(c+1)p\big)\frac{(n-s-cp)^{2}}{\frac{n-s-(c-6)p}{2p}-2\sqrt{\frac{n-s-(c-4)p}{2p}}}<2(c-1)(n-\frac{3}{2}s-2p)\frac{(n-s-cp)^{2}}{\frac{n-s-cp}{2p}}\nonumber
\\\Leftrightarrow& \frac{n-2s-2(c+1)p}{n-s-(c-6)p-2\sqrt{2p}\sqrt{n-s-(c-4)p}}<\frac{2(c-1)\big(n-2s-2(c+1)p+\frac{1}{2}s+2cp\big)}{n-s-cp}\nonumber
\\\Leftrightarrow& \frac{n-s-cp}{n-s-(c-6)p-2\sqrt{2p}\sqrt{n-s-(c-4)p}}<\frac{2(c-1)\big(n-2s-2(c+1)p+\frac{1}{2}s+2cp\big)}{n-2s-2(c+1)p}\nonumber
\\\Leftrightarrow& \frac{1}{1+\frac{6p}{n-s-cp}-2\sqrt{2p}\sqrt{\frac{1}{n-s-cp}+\frac{4p}{(n-s-cp)^{2}}}}<\frac{\frac{1}{2}s+2cp}{n-2s-2(c+1)p}+2(c-1).\nonumber
\end{align}

Since $\frac{1}{1+\frac{6p}{n-s-cp}-2\sqrt{2p}\sqrt{\frac{1}{n-s-cp}+\frac{4p}{(n-s-cp)^{2}}}}$ attains its maximum when $n-s-cp=4p$, we get $\frac{1}{1+\frac{6p}{4p}-2\sqrt{2p}\sqrt{\frac{1}{4p}+\frac{4p}{(4p)^{2}}}}=2\leq2(c-1)<\frac{\frac{1}{2}s+2cp}{n-2s-2(c+1)p}+2(c-1)$. Thus, the inequality (1) holds.

The inequality (2) implies
\begin{align} &r^{2}-(n-s-cp)s-4(c-1)p(n-s-cp)<0\nonumber
\\\Leftrightarrow&
\frac{(n-s-cp)^{2}}{\frac{n-s-(c-6)p}{2p}-2\sqrt{\frac{n-s-(c-4)p}{2p}}}<\big(2(c-1)+\frac{s}{2p}\big)\frac{(n-s-cp)^{2}}{\frac{n-s-cp}{2p}}\nonumber
\\\Leftrightarrow&
\frac{n-s-cp}{n-s-(c-6)p-2\sqrt{2p}\sqrt{n-s-(c-4)p}}<2(c-1)+\frac{s}{2p}\nonumber
\\\Leftrightarrow& \frac{1}{1+\frac{6p}{n-s-cp}-2\sqrt{2p}\sqrt{\frac{1}{n-s-cp}+\frac{4p}{(n-s-cp)^{2}}}}<2(c-1)+\frac{s}{2p}.\nonumber
\end{align}
Since $\frac{1}{1+\frac{6p}{4p}-2\sqrt{2p}\sqrt{\frac{1}{4p}+\frac{4p}{(4p)^{2}}}}=2<2(c-1)+\frac{s}{2p}$, we obtain the inequality (2) holds.

Thus, $det\big((2n-s-2p-2+r)I-B\big)<0$ and $\eta_{1}(G_{1})=\rho(B)>2n-s-2p-2+r$, which implies $\eta_{1}(G_{2})-\eta_{1}(G_{1})>0$.

If $(n_{1},n_{2},\ldots,n_{c})=(n-s-(c-1)p,p,\ldots,p)$, then $K_{s}\vee (K_{n_{1}}+K_{n_{2}}+\cdots+K_{n_{c}})\cong K_{s}\vee (K_{n-s-(c-1)p}+(c-1)K_{p})$, and $\eta_{1}(G_{1})=\eta_{1}(G_{2})$.
Conversely, if $\eta_{1}(G_{1})=\eta_{1}(G_{2})$, then $n_{2}=n_{3}=\cdots=n_{c}=p$, and $(n_{1},n_{2},\ldots,n_{c})=(n-s-(c-1)p,p,\ldots,p)$.

This completes the proof.
$\hfill\square$\\

Now, we can give a proof of Theorem \ref{th:1.1.}.\\

\noindent\textbf{Proof of Theorem \ref{th:1.1.}.}

Suppose to the contrary that $G$ is not 1-tough.
Thus, there exists a vertex set $S\subseteq V(G)$ such that $c(G-S)>|S|$.
Let $|S|=s$ and $c(G-S)=c$, then $c\geq s+1$ and $n\geq2s+1$.

We know $G$ is a spanning subgraph of $G^{\prime}=K_{s}\vee(K_{n_{1}}+K_{n_{2}}+\cdots+K_{n_{s+1}})$ for some integers $n_{1}\geq n_{2}\geq\cdots\geq n_{s+1}\geq1$ and $\sum\limits_{i=1}^{s+1}n_{i}=n-s$.
By Lemma \ref{le:2.1.}, we obtain
\begin{align}
\eta_{1}(G^{\prime})\leq\eta_{1}(G)
\end{align}
with equality holds  if and only if $G\cong G^{\prime}$.
We divide the following proof into three cases.\\

\textbf{Case 1.} $s\geq\delta+1$.

Let $G_{1}=K_{s}\vee(K_{n-2s}+sK_{1})$.
By Lemma \ref{le:2.4.},
\begin{align}
\eta_{1}(G_{1})\leq\eta_{1}(G^{\prime})
\end{align}
with equality holds  if and only if $G_{1}\cong G^{\prime}$.

Let $G^{\ast}=K_{\delta}\vee(K_{n-2\delta}+\delta K_{1})$.
The distance signless Laplacian matrix of $G^{\ast}$ is $Q_{D}(G^{\ast})=D(G^{\ast})+Tr(G^{\ast})$, where
\begin{align*}
D(G^{\ast})=
\bordermatrix{%
& \delta & n-2\delta & \delta\cr
\delta & 2(J-I) & 2J & J  \cr
n-2\delta & 2J & J-I & J  \cr
\delta & J & J & J-I  \cr
},
\end{align*}
\begin{align*}
Tr(G^{\ast})=
\bordermatrix{%
& \delta & n-2\delta & \delta\cr
\delta & (2n-\delta-2)I & 0 & 0  \cr
n-2\delta & 0 & (n+\delta-1)I & 0  \cr
\delta & 0 & 0 & (n-1)I  \cr
}.
\end{align*}

The quotient matrix of the distance signless Laplacian matrix $Q_{D}(G^{\ast})$ according to the partition $\{V(\delta K_{1}), V(K_{n-2\delta}), V(K_{\delta})\}$ is
\begin{align*}
B^{\ast}=
\left(\begin{matrix}
 2n+\delta-4 & 2n-4\delta & \delta  \cr
 2\delta & 2n-\delta-2 & \delta  \cr
 \delta & n-2\delta & n+\delta-2  \cr
\end{matrix}\right),
\end{align*}
where the partition is equitable.

According to the definition of the Wiener index,
\begin{align*} W(G^{\ast})&=\sum\limits_{i<j}d_{ij}(G^{\ast})\nonumber
\\&=\frac{2[1+(\delta-1)](\delta-1)}{2}+2(n-2\delta)\delta+\delta^{2}+\frac{[1+(n-\delta-1)](n-\delta-1)}{2}\nonumber
\\&=\frac{1}{2}n^{2}+(\delta-\frac{1}{2})-\frac{3}{2}\delta^{2}-\frac{1}{2}\delta.\nonumber
\end{align*}
By  Lemma \ref{le:2.3.}, since $n\geq11\delta$ and $\delta\geq2$,
\begin{align*} \eta_{1}(G^{\ast})=\frac{4W(G^{\ast})}{n}&=\frac{2n^{2}+(4\delta-2)n-6\delta^{2}-2\delta}{n}\nonumber
\\&=2n+4\delta-2-\frac{6\delta^{2}+2\delta}{n}\nonumber
\\&\geq(2n+3\delta-2)+\delta-\frac{6\delta^{2}+2\delta}{11\delta}\nonumber
\\&=(2n+3\delta-2)+\frac{5\delta^{2}-2\delta}{11\delta}\nonumber
\\&> 2n+3\delta-2.\nonumber
\end{align*}

The characteristic polynomial of $B^{\ast}$ is
\begin{align*}
P(B^{\ast},x)=&x^{3}-(5n+\delta-8)x^{2}+\big(8n^{2}-(\delta+26)n+8\delta^{2}-4\delta+20\big)x\nonumber
\\&-4n^{3}+(2\delta+20)n^{2}-(2\delta^{2}+2\delta+32)n-2\delta^{3}+18\delta^{2}-4\delta+16.\nonumber
\end{align*}

By  Lemma \ref{le:2.2.}, $\eta_{1}(G^{\ast})=\rho(B^{\ast})$ is the largest root of the equation $P(B^{\ast},x)=0$.
$Q_{D}(G_{1})$ has the equitable matrix $B_{1}$ which is obtained by replacing $\delta$ with $s$ in $B^{\ast}$, and $\eta_{1}(G_{1})=\rho(B_{1})$ is the largest root of the equation $P(B_{1},x)=0$.
Then
\begin{align*}
P(B^{\ast},x)&-P(B_{1},x)=(s-\delta)\Big[x^{2}+(n-8s-8\delta+4)x\nonumber
\\&+2s^{2}+8ns+2\delta s-18s-2n^{2}+8\delta n+2n+2\delta^{2}-18\delta+4\Big].\nonumber
\end{align*}

Let $f(x)=x^{2}+(n-8s-8\delta+4)x
+2s^{2}+8ns+2\delta s-18s-2n^{2}+8\delta n+2n+2\delta^{2}-18\delta+4$, and $P(B^{\ast},x)-P(B_{1},x)=(s-\delta)f(x)$.
Since $s\geq\delta+1$ and $n\geq2s+1$, then $\delta+1\leq s\leq\frac{n-1}{2}$.
The symmetry axis of $f(x)$ is
\begin{align*}
x_{0}&=-\frac{n-8s-8\delta+4}{2}\nonumber
\\&=(2n-\frac{5}{2}n)+4s+4\delta-2\nonumber
\\&=(2n+3\delta-2)-\frac{5}{2}n+4s+\delta\nonumber
\\&\leq(2n+3\delta-2)-\frac{5}{2}(2s+1)+4s+\delta\nonumber
\\&=(2n+3\delta-2)-s+\delta-\frac{5}{2}\nonumber
\\&\leq(2n+3\delta-2)-\frac{7}{2}\nonumber
\\&<2n+3\delta-2.\nonumber
\end{align*}

Thus, $f(x)$ is monotonically increasing with respect to $x\in\left[2n+3\delta-2, +\infty\right)$.
Since $\delta+1\leq s\leq\frac{n-1}{2}$, $n\geq11\delta$ and $\delta\geq2$, we obtain
\begin{align*}
f(x)\geq& f(2n+3\delta-2)\nonumber
\\=&2s^{2}-(8n+22\delta+2)s+4n^{2}+7\delta n-13\delta^{2}-2\delta\nonumber
\\\geq&2(\frac{n-1}{2})^{2}-(8n+22\delta+2)\frac{n-1}{2}+4n^{2}+7\delta n-13\delta^{2}-2\delta\nonumber
\\=&\frac{1}{2}n^{2}-(4\delta-2)n-13\delta^{2}+9\delta+\frac{3}{2}\nonumber
\\\geq&\frac{7}{2}\delta^{2}+31\delta+\frac{3}{2}>0.\nonumber
\end{align*}
Thus, $P(B^{\ast},x)>P(B_{1},x)$ for $x\in\left[2n+3\delta-2, +\infty\right)$.

Combining with $\eta_{1}(G^{\ast})> 2n+3\delta-2$, we have
\begin{align}
\eta_{1}(G^{\ast})<\eta_{1}(G_{1}).
\end{align}

By inequalities (3), (4) and (5), we can see that $\eta_{1}(G^{\ast})<\eta_{1}(G_{1})\leq\eta_{1}(G^{\prime})\leq\eta_{1}(G)$, which is a contradiction.\\

\textbf{Case 2.} $s=\delta$.

Let $G^{\prime}=K_{\delta}\vee(K_{n_{1}}+K_{n_{2}}+\cdots+K_{n_{\delta+1}})$ and $G^{\ast}=K_{\delta}\vee(K_{n-2\delta}+\delta K_{1})$.
By   Lemma \ref{le:2.4.}, we obtain
\begin{align}
\eta_{1}(G^{\ast})\leq\eta_{1}(G^{\prime})
\end{align}
with equality holds if and only if $G^{\prime}\cong G^{\ast}$.
By equalities (3) and (6), we obtain
$$\eta_{1}(G^{\ast})\leq\eta_{1}(G),$$ with equality holds if and only if  $G\cong G^{\ast}$.

Due to the assumption $\eta_{1}(G^{\ast})\geq\eta_{1}(G)$, we can get $\eta_{1}(G^{\ast})=\eta_{1}(G)$ which implies $G\cong G^{\ast}$.
Let $S=V(K_{\delta})$.
Thus, $$\frac{s}{c}=\frac{\delta}{\delta+1}<1,$$ and then $\tau(G^{\ast})<1$, which implies $G^{\ast}=K_{\delta}\vee(K_{n-2\delta}+\delta K_{1})$ is not 1-tough.
So $G\cong G^{\ast}$.\\

\textbf{Case 3.} $1\leq s\leq\delta-1$.

Let $G_{2}=K_{s}\vee(K_{n-s-s(\delta-s+1)}+sK_{\delta-s+1})$.
Recall that $G$ is a spanning subgraph of $G^{\prime}=K_{s}\vee(K_{n_{1}}+K_{n_{2}}+\cdots+K_{n_{s+1}})$, where $n_{1}\geq n_{2}\geq\cdots\geq n_{s+1}\geq1$ and $\sum\limits_{i=1}^{s+1}n_{i}=n-s$.
Since the minimum degree of $G^{\prime}$ is at least $\delta$, then $n_{s+1}\geq\delta-s+1$.
By Lemma \ref{le:2.5.}, we have
\begin{align}
\eta_{1}(G_{2})\leq\eta_{1}(G^{\prime})
\end{align}
with equality holds if and only if $(n_{1}, n_{2},\ldots, n_{s+1})=(n-s-s(\delta-s+1),\delta-s+1,\ldots,\delta-s+1)$.\\

\textbf{Subcase 3.1.} $s=1$.

In this case, $G_{2}=K_{1}\vee(K_{n-\delta-1}+K_{\delta})$ and the distance signless Laplacian matrix of $G_{2}$ is $Q_{D}(G_{2})=D(G_{2})+Tr(G_{2})$, where
\begin{align*}
D(G_{2})=
\bordermatrix{%
& \delta & n-\delta-1 & 1\cr
\delta & J-I & 2J & J  \cr
n-\delta-1 & 2J & J-I & J  \cr
1 & J & J & 0  \cr
},
\end{align*}
\begin{align*}
Tr(G_{2})=
\bordermatrix{%
& \delta & n-\delta-1 & 1\cr
\delta & (2n-\delta-2)I & 0 & 0  \cr
n-\delta-1 & 0 & (n+\delta-1)I & 0  \cr
1 & 0 & 0 & (n-1)I  \cr
}.
\end{align*}

Recall that $G^{\ast}=K_{\delta}\vee(K_{n-2\delta}+\delta K_{1})$.
Let $\textbf{x}$ be the Perron vector of $Q_{D}(G^{\ast})$ corresponding to $\eta_{1}(G^{\ast})$, and we take $x_{u}=x_{1}$ for all $u\in V(\delta K_{1})$, $x_{v}=x_{2}$ for all $x_{v}\in V(K_{n-2\delta})$, and $x_{w}=x_{3}$ for all $x_{w}\in V( K_{\delta})$.
Due to $Q_{D}(G^{\ast})\textbf{x}=\eta_{1}(G^{\ast})\textbf{x}$,
$$\left\{
\begin{aligned}
&\eta_{1}(G^{\ast})x_{1}=(2n+\delta-4)x_{1}+(2n-4\delta)x_{2}+\delta x_{3}\\
&\eta_{1}(G^{\ast})x_{3}=\delta x_{1}+(n-2\delta)x_{2}+(n+\delta-2)px_{3}.
\end{aligned}
\right.$$
Thus,
\begin{align*}
x_{3}=\frac{\eta_{1}(G^{\ast})-2n+\delta+4}{2\eta_{1}(G^{\ast})-2n-\delta+4}x_{1}=(\frac{1}{2}+\frac{-n+\frac{3}{2}\delta+2}{2\eta_{1}(G^{\ast})-2n-\delta+4})x_{1}.
\end{align*}
Let $k=\frac{-n+\frac{3}{2}\delta+2}{2\eta_{1}(G^{\ast})-2n-\delta+4}$, we have $x_{3}=(\frac{1}{2}+k)x_{1}$ and $k=\frac{-n+\frac{3}{2}\delta+2}{2\eta_{1}(G^{\ast})-2n-\delta+4}<\frac{-n+\frac{3}{2}\delta+2}{2n+5\delta}=-\frac{1}{2}+\frac{4\delta+2}{2n+5\delta}<0$, since $\eta_{1}(G^{\ast})>2n+3\delta-2$ and $n\geq11\delta$.

Let $\textbf{\~{x}}=(\underbrace{x_{1}, \ldots, x_{1}}_{\delta}, \underbrace{x_{2}, \ldots, x_{2}}_{n-\delta}, \underbrace{x_{3}, \ldots, x_{3}}_{\delta})^{T}$.
The matrix $Q_{D}(G_{2})-Q_{D}(G^{\ast})$ is as follows:

\begin{align*}
\bordermatrix{%
& \delta & n-2\delta & \delta-1 & 1 \cr
\delta     & -(J-I) & 0 & J & 0  \cr
n-2\delta & 0 & 0 & 0 & 0 \cr
\delta-1  & J & 0 & \delta &0 \cr
1& 0 & 0 & 0 & 0\cr
}.
\end{align*}

Thus,
\begin{align*}
\eta_{1}(G_{2})-\eta_{1}(G^{\ast})&\geq \textbf{\~{x}}^{T}(Q_{D}(G_{2})-Q_{D}(G^{\ast}))\textbf{\~{x}}
\\&=-\delta(\delta-1)x_{1}^{2}+2\delta(\delta-1)x_{1}x_{3}+\delta(\delta-1) x_{3}^{2}
\\&=\delta(\delta-1)[x_{3}^{2}+(2x_{3}-x_{1})x_{1}]
\\&=\delta(\delta-1)(x_{3}^{2}+2kx_{1}^{2})
\\&=(k^{2}+3k+\frac{1}{4})x_{1}^{2}
\\&>0.
\end{align*}

Thus, we obtain
\begin{align}
\eta_{1}(G^{\ast})<\eta_{1}(G_{2}).
\end{align}\\

\textbf{Subcase 3.2.} $2\leq s\leq\delta-1$.

The distance signless Laplacian matrix of $G_{2}$ is $Q_{D}(G_{2})=D(G_{2})+Tr(G_{2})$, where
\begin{align*}
\scalebox{0.8}{$
D(G_{2})=
\bordermatrix{%
& \delta-s+1 & \cdots & \delta-s+1 & n-s-s(\delta-s+1)  & s\cr
\delta-s+1 & J-I & \cdots & 2J & 2J & J \cr
\vdots & \vdots & \ddots & \vdots & \vdots & \vdots \cr
\delta-s+1 & 2J & \cdots & J-I & 2J & J \cr
n-s-s(\delta-s+1) & 2J & \cdots & 2J & J-I & J \cr
s & J & \cdots & J & J & J-I \cr
},$}
\end{align*}
\begin{align*}
\scalebox{0.8}{$
Tr(G_{2})=
\bordermatrix{%
& \delta-s+1 & \cdots & \delta-s+1 & n-s-s(\delta-s+1)  & s\cr
\delta-s+1 & (2n-\delta-2)I & \cdots & 0 & 0 & 0 \cr
\vdots & \vdots & \ddots & \vdots & \vdots & \vdots \cr
\delta-s+1 & 0 & 0 & (2n-\delta-2)I & 0 & 0 \cr
n-s-s(\delta-s+1) & 0 & 0 & 0 & (n+s(\delta-s+1)-1)I &0 \cr
s & 0 & 0 & 0 & 0 & (n-1)I \cr
}.$}
\end{align*}

The quotient matrix of the distance signless Laplacian matrix $Q_{D}(G_{2})$ according to the partition $\{V(s K_{\delta-s+1}), V(K_{n-s-s(\delta-s+1)}), V(K_{s})\}$ of $G_{2}$ is
\begin{align*}
\scalebox{0.9}{$
B_{2}=
\left(\begin{matrix}
 2n-s+2(s-1)(\delta-s+1)-2 & 2[n-s-s(\delta-s+1)] & s  \cr
 2s(\delta-s+1) & 2n-s-2 & s  \cr
 s(\delta-s+1) & n-s-s(\delta-s+1) & n+s-2  \cr
\end{matrix}\right).$}
\end{align*}
The partition is equitable and we can obtain the characteristic polynomial of $B_{2}$,
{\small\begin{align*}
P(B_{2},x)=&x^{3}+\big(2s^{2}-(2\delta+3)s-5n+2\delta+8\big)x^{2}+\big(4s^{2}-(8\delta+12)s^{3}
+(-2n+4\delta^{2}+12\delta+16)s^{2}\nonumber
\\&+(2\delta n+5n-8\delta-12)s+8n^{2}-(6\delta+26)n+8\delta+20\big)x-2s^{5}-(4n-4\delta-14)s^{4}\nonumber
\\&+(8\delta n+12n-2\delta^{2}-24\delta-30)s^{3}-(4\delta^{2}n+12\delta n+12n-10\delta^{2}-28\delta-26)s^{2}\nonumber
\\&-(2n^{2}-4\delta n-10n+8\delta+12)s-4n^{3}+4\delta n^{2}+20n^{2}-12\delta n-32n+8\delta+16.\nonumber
\end{align*}}

By Lemma \ref{le:2.2.}, $\eta_{1}(G_{2})=\rho(B_{2})$ is the largest root of the equation $P(B_{2},x)=0$.
Since $\eta_{1}(G^{\ast})>2n+3\delta-2$,
\begin{align*}
&P(B^{\ast},2n+3\delta-2)-P(B_{2},2n+3\delta-2)
\\=&(\delta-s)\Big[4(s-1)n^{2}+\big(4s^{3}-(4\delta+12)s^{2}+(18\delta+8)s-13\delta\big)n
\\&-2s^{4}+(14\delta+6)s^{3}-(12\delta^{2}+38\delta-6)s^{2}+(18\delta^{2}+22\delta+2)s-5\delta^{2}-2\delta\Big]
\\=&(\delta-s)g(n).
\end{align*}

Note that $2\leq s\leq\delta-1$ and $\delta\geq s+1\geq3$.
Thus, the symmetry axis of $g(n)$ is
\begin{align*}
n_{0}&=-\frac{4s^{3}-(4\delta+12)s^{2}+(18\delta+8)s-13\delta}{8(s-1)}\nonumber
\\&=-\frac{4(s-1)^{3}-4\delta(s-1)^{2}+(10\delta-4)(s-1)+3\delta}{8(s-1)}\nonumber
\\&=-\frac{1}{4}\big[2(s-1)^{2}-2\delta(s-1)+(5\delta-2)\big]-\frac{3\delta}{8(s-1)}\nonumber
\\&<-\frac{1}{2}(s-1)^{2}+\frac{\delta}{2}(s-1)-\frac{5\delta-2}{4}.\nonumber
\end{align*}

If $3\leq\delta\leq4$, then
\begin{align*}
n_{0}&<-\frac{1}{2}(s-1)^{2}+\frac{\delta}{2}(s-1)-\frac{5\delta-2}{4}\nonumber
\\&\leq-\frac{1}{2}(\delta-2)^{2}+\frac{\delta}{2}(\delta-2)-\frac{5\delta-2}{4}\nonumber
\\&=-\frac{1}{4}\delta-\frac{3}{2}<\frac{1}{2}\delta^{2}+2\delta.\nonumber
\end{align*}

If $\delta\geq5$, then
\begin{align*}
n_{0}&<-\frac{1}{2}(s-1)^{2}+\frac{\delta}{2}(s-1)-\frac{5\delta-2}{4}\nonumber
\\&\leq-\frac{1}{2}(\frac{\delta}{2})^{2}+(\frac{\delta}{2})^{2}-\frac{5\delta-2}{4}\nonumber
\\&=\frac{1}{8}\delta^{2}-\frac{5\delta-2}{4}<\frac{1}{2}\delta^{2}+2\delta.\nonumber
\end{align*}

This implies that $g(n)$ is monotonically increasing with respect to $n\in\left[\frac{1}{2}\delta^{2}+2\delta, +\infty\right)$.

Thus,
\begin{align*}
g(n)\geq& g(\frac{1}{2}\delta^{2}+2\delta)
\\=&\delta\Big[\delta\big((s-1)\delta^{2}-(2s^{2}-2s-3)\delta+2s^{3}-6s^{2}+4s-2\big)+2s^{3}-2s^{2}-2s+2\Big]\nonumber
\\&-2s^{4}+6s^{3}-6s^{2}+2s
\\\geq&\delta\Big[\delta\big(s^{3}-5s^{2}+8s\big)+2s^{3}-2s^{2}-2s+2\Big]-2s^{4}+6s^{3}-6s^{2}+2s
\\\geq&\delta\Big[\big(s^{4}-4s^{3}+3s^{2}+8s\big)+2s^{3}-6s^{2}+4s-2\Big]-2s^{4}+6s^{3}-6s^{2}+2s
\\\geq&s^{5}-s^{4}+7s^{3}+21s^{2}+12s-2s^{4}+6s^{3}-6s^{2}+2s
\\=&s^{5}-3s^{4}+13s^{3}+15s^{2}+14s
\\>&0.
\end{align*}

By the above analyses, we can see that \begin{align}
P(B^{\ast},2n+3\delta-2)-P(B_{2},2n+3\delta-2)>0.
\end{align}

Since $\eta_{1}(G^{\ast})>2n+3\delta-2$, then $x\in \left[2n+3\delta-2,+\infty\right)$, and $s\geq2$.
Take the derivative of the characteristic polynomials and we have
\begin{align*}
&P^{\prime}(B^{\ast},x)-P^{\prime}(B_{2},x)
\\&=(\delta-s)\Big[(4s-6)x+4s^{3}-(4\delta+12)s^{2}-(2n-16)s+5n+8\delta-12\Big]
\\&\geq(\delta-s)\Big[(4s-6)(2n+3\delta-2)+4s^{3}-(4\delta+12)s^{2}-(2n-16)s+5n+8\delta-12\Big]
\\&=(\delta-s)\Big[4s^{3}-(4\delta+12)s^{2}+(6n+12\delta+8)s-7n-10\delta\Big]
\\&=(\delta-s)p(s).
\end{align*}

Thus $p(s)=4s^{3}-(4\delta+12)s^{2}+(6n+12\delta+8)s-7n-10\delta$ and $2\leq s\leq\delta-1$.
We obtain that $p^{\prime}(s)=12s^{2}-8(\delta+3)s+6n+12\delta+8$ and the symmetry axis of $p^{\prime}(s)$ is $s_{0}=\frac{\delta+3}{3}$.
Since $n\geq\frac{1}{2}\delta^{2}+2\delta$ and $\delta\geq3$, we have $$p^{\prime}(s)\geq p^{\prime}(\frac{\delta+3}{3})\geq 3\delta^{2}+12\delta-\frac{4}{3}\delta^{2}+4\delta-4=\frac{5}{3}\delta^{2}+16\delta-4>0,$$
which implies $p(s)$ is monotonically increasing for $2\leq s\leq\delta-1$.
Combining with $n\geq11\delta$ and $\delta\geq3$, we obtain that
$p(s)\geq p(2)=5n-2\delta>0$.
Thus,
\begin{align}
P^{\prime}(B^{\ast},x)-P^{\prime}(B_{2},x)>0.
\end{align}

Moreover, since $P^{\prime}(B^{\ast},x)=3x^{2}-2(5n+\delta-8)x+8n^{2}-(\delta+26)n+8\delta^{2}-4\delta+20$,
the symmetry axis of $P^{\prime}(B^{\ast},x)$ is \begin{align*}
x_{0}=\frac{5n+\delta-8}{3}&=(2n+3\delta-2)-\frac{1}{3}n-\frac{8}{3}\delta-\frac{2}{3}
\\&\leq(2n+3\delta-2)-\frac{1}{3}(11\delta)-\frac{8}{3}\delta-\frac{2}{3}
\\&\leq(2n+3\delta-2)-\frac{19}{3}\delta-\frac{2}{3}
\\&<2n+3\delta-2.
\end{align*}
Thus, we have $P^{\prime}(B^{\ast},x)\geq P^{\prime}(B^{\ast},2n+3\delta-2)=(\delta+2)n+29\delta^{2}+12\delta>0$.
It follows that $P(B^{\ast},x)$ is monotonically increasing with respect to $x\in \left[2n+3\delta-2,+\infty\right)$.
Combining  with inequalities (9) and (10), we deduce that
\begin{align}
\eta_{1}(G^{\ast})<\eta_{1}(G_{2}).
\end{align}
By inequalities (3), (7), (8) and (11), we can see that $\eta_{1}(G^{\ast})<\eta_{1}(G_{2})\leq\eta_{1}(G^{\prime})\leq\eta_{1}(G)$, which is a contradiction.

By the above analyses, we can see that $\eta_{1}(G^{\ast})<\eta_{1}(G)$, a contradiction.

This completes the proof.
$\hfill\square$\\

\section{Proof of Theorem \ref{th:1.2.} }
\hspace{1.3em}

In this section, we give the proof of Theorem \ref{th:1.2.} and the following lemma which will be used further.

\noindent\begin{lemma}\label{le:3.1.}\cite{AH}
If $G$ is a connected graph of order $n\geq2$, then
$$\eta_{1}(G)\geq\eta_{1}(K_{n})=2n-2$$
with equality if and only if $G$ is the complete graph $K_{n}$.
\end{lemma}

Now, we present sufficient condition for graphs to be $t$-tough, where $\frac{1}{t}$ is a positive integer.\\

\noindent\textbf{Proof of Theorem\ref{th:1.2.}.}

Assume to the contrary that $G$ is not $t$-tough.
There exists a vertex subset $S\subseteq V(G)$ such that $tc(G-S)>|S|$.
Let $|S|=s$ and $c(G-S)=c$, then $tc>s$.

When $\frac{1}{t}$ is a positive number, we obtain $c\geq\frac{s}{t}+1$.
We know $G$ is a  spanning subgraph of $G^{\prime}=K_{s}\vee(K_{n_{1}}+K_{n_{2}}+\cdots+K_{n_{\frac{s}{t}+1}})$ for some integers $n_{1}\geq n_{2}\geq\cdots\geq n_{\frac{s}{t}+1}\geq1$ and $\sum\limits_{i=1}^{\frac{s}{t}+1}n_{i}=n-s$.
By Lemma \ref{le:2.1.}, we obtain
\begin{align}
\eta_{1}(G^{\prime})\leq\eta_{1}(G)
\end{align}
with equality holds  if and only if $G\cong G^{\prime}$.

Let $\widetilde{G}=K_{s}\vee(K_{n-\frac{s}{t}-s}+\frac{s}{t}K_{1})$.
 According to  Lemma \ref{le:2.4.},  we have
\begin{align}
\eta_{1}(\widetilde{G})\leq\eta_{1}(G^{\prime})
\end{align}
with equality holds  if and only if $\widetilde{G}\cong G^{\prime}$.

Next, we consider the following two cases based on the value of $s$.\\

\textbf{Case 1.} $s=1$.

In this case, $\widetilde{G}=K_{1}\vee(K_{n-\frac{1}{t}-1}+\frac{1}{t}K_{1})$.
By inequalities (12) and (13), we obtain $\eta_{1}(K_{1}\vee(K_{n-\frac{1}{t}-1}+\frac{1}{t}K_{1}))\leq\eta_{1}(G)$ with equality holding if and only if $G\cong K_{1}\vee(K_{n-\frac{1}{t}-1}+\frac{1}{t}K_{1})$.
Since the assumption $\eta_{1}(G)\leq\eta_{1}(K_{1}\vee(K_{n-\frac{1}{t}-1}+\frac{1}{t}K_{1}))$, we have $\eta_{1}(G)=\eta_{1}(K_{1}\vee(K_{n-\frac{1}{t}-1}+\frac{1}{t}K_{1}))$ which implies $G\cong K_{1}\vee(K_{n-\frac{1}{t}-1}+\frac{1}{t}K_{1})$.

Take $S=V(K_{1})$, thus $\tau(K_{1}\vee(K_{n-\frac{1}{t}-1}+\frac{1}{t}K_{1}))=\frac{s}{c}=\frac{1}{1+\frac{1}{t}}<t$, which implies $K_{1}\vee(K_{n-\frac{1}{t}-1}+\frac{1}{t}K_{1})$ is not $t$-tough.
So $G\cong K_{1}\vee(K_{n-\frac{1}{t}-1}+\frac{1}{t}K_{1})$.\\

\textbf{Case 2.} $s\geq2$.

According to $\widetilde{G}=K_{s}\vee(K_{n-\frac{s}{t}-s}+\frac{s}{t}K_{1})$, the equitable quotient matrix of $Q(\widetilde{G})$ is as follows:
\begin{align*}
\widetilde{B}=
\left(\begin{matrix}
 n+s-2 & n-s-\frac{s}{t} & \frac{s}{t}  \cr
 s & 2n-s-2 & \frac{2s}{t}  \cr
 s & 2n-2s-\frac{2s}{t} & 2n-s+\frac{2s}{t}-4  \cr
\end{matrix}\right).
\end{align*}

Thus, the characteristic polynomial of $\widetilde{B}$ is
\begin{align}
P(\widetilde{B},x)=&x^{3}-(5n+\frac{2}{t}s-s-8)x^{2}+\big((\frac{4}{t^{2}}+\frac{4}{t})s^{2}+(\frac{2n}{t}-3n-\frac{8}{t}+4)s+8n^{2}-26n+20\big)x\nonumber
\\&-4n^{3}+2(s+10)n^{2}-(\frac{4}{t^{2}}s^{2}+\frac{4}{t}s^{2}-\frac{4}{t}s+6s+32)n\nonumber
\\&-\frac{2}{t^{2}}s^{3}+\frac{8}{t^{2}}s^{2}+\frac{10}{t}s^{2}-\frac{8}{t}s+4s+16.\nonumber
\end{align}

Let $G^{\ast}=K_{1}\vee(K_{n-\frac{1}{t}-1}+\frac{1}{t}K_{1})$.
Note that $Q_{D}(G^{\ast})$ has the equitable quotient matrix $B_{t}$, which is obtained by taking $s=1$ in $\widetilde{B}$.
Hence,
\begin{align*}
P(B_{t},x)-P(\widetilde{B},x)=&(s-1)\Big[(\frac{2}{t}-1)x^{2}-\big(\frac{4}{t^{2}}(s+1)+\frac{2}{t}(n+2s-2)-3n+4\big)x
\\&+\frac{2}{t^{2}}\big(s^{2}+(2n-3)(s+1)\big)+\frac{2}{t}\big((2n-5)s-1\big)-2n^{2}+6n-4\Big]
\\=&(s-1)h(x).\end{align*}

According to Lemma \ref{le:3.1.} and $n\geq s+\frac{s}{t}+1$, then $2\leq s\leq\frac{n-1}{\frac{1}{t}+1}$.
The symmetry axis of $h(x)$ is
\begin{align*}
x_{0}=&\frac{\frac{4}{t^{2}}(s+1)+\frac{2}{t}(n+2s-2)-3n+4}{\frac{4}{t}-2}
\\\leq&\frac{\frac{4}{t}(n-1)+\frac{4}{t^{2}}+(\frac{2}{t}-3)n-\frac{4}{t}+4}{\frac{4}{t}-2}
\\=&\frac{3}{2}n+\frac{1}{\frac{4}{t}-2}+\frac{1}{t}-\frac{3}{2}
\\\leq&2n-\frac{1}{2}s-\frac{1}{2t}s-\frac{1}{2}+\frac{1}{\frac{4}{t}-2}+\frac{1}{t}-\frac{3}{2}
\\\leq&2n-1-\frac{1}{t}-\frac{1}{2}+\frac{1}{\frac{4}{t}-2}+\frac{1}{t}-\frac{3}{2}
\\\leq&2n+\frac{1}{\frac{4}{t}-2}-3<2n-2.
\end{align*}
This implies that $h(x)$ is monotonically increasing with respect to $x\in \left[2n-2,+\infty\right)$.

Since $n\geq \frac{2}{t^{2}}+\frac{3}{t}+3t+4$ and $2\leq s\leq\frac{n-1}{\frac{1}{t}+1}$, then
\begin{align*}
h(x)\geq& h(2n-2)
\\=&\frac{2}{t^{2}}s^{2}-\big((\frac{4}{t^{2}}+\frac{4}{t})n-\frac{2}{t^{2}}+\frac{2}{t}\big)s-\frac{2}{t^{2}}(2n-1)-\frac{2}{t}(2n+1)+2n^{2}
\\\geq&\frac{2}{t^{2}}(\frac{n-1}{\frac{1}{t}+1})^{2}-\big((\frac{4}{t^{2}}+\frac{4}{t})n-\frac{2}{t^{2}}+\frac{2}{t}\big)\frac{n-1}{\frac{1}{t}+1}-\frac{2}{t^{2}}(2n-1)-\frac{2}{t}(2n+1)+2n^{2}
\\=&\frac{1}{t^{2}(1+t)^{2}}\Big[2t^{2}n^{2}-2\big(t^{3}+4t^{2}+3t+2\big)n+2\Big]
\\\geq&\frac{1}{t^{2}(1+t)^{2}}\Big[2t^{2}\big(\frac{2}{t^{2}}+\frac{3}{t}+3t+4\big)^{2}-2\big(t^{3}+4t^{2}+3t+2\big)\big(\frac{2}{t^{2}}+\frac{3}{t}+3t+4\big)+2\Big]
\\=&\frac{12t^{4}+16t^{3}+12t^{2}+8t+2}{t^{2}(1+t)^{2}}
\\>&0.
\end{align*}

Combining $h(x)>0$ and $s\geq2$, we deduce that $P(B_{t},x)-P(\widetilde{B},x)>0$ for $x\in \left[2n-2,+\infty\right)$ which implies $\eta_{1}(B_{t})>\eta_{1}(\widetilde{B})$.
Thus, combining with Lemma \ref{le:2.2.}, we have
\begin{align}
\eta_{1}(G^{\ast})\leq\eta_{1}(\widetilde{G}).
\end{align}

By inequalities (12), (13) and (14), we can see that $\eta_{1}(G^{\ast})<\eta_{1}(\widetilde{G})\leq\eta_{1}(G^{\prime})\leq\eta_{1}(G)$, which is a contradiction.

By the above analyses, we can see that $\eta_{1}(G^{\ast})<\eta_{1}(G)$, a contradiction.

This completes the proof.
$\hfill\square$\\

\section{Proof of Theorem \ref{th:1.3.} }
\hspace{1.3em}

Chen et al. \cite{CLX} considered the new parameter $\tau'(G)=\min\{\frac{|S|}{c(G-S)-1}\}$ and established two lower bounds on the size to guarantee that a graph $G$ is $t$-tough.
In order to proof Theorem \ref{th:1.3.}, we list their result.

\noindent\begin{lemma}\label{le:4.1.}\cite{CLX}
Let $G$ be a connected graph of order $n$ and minimum degree $\delta$, then the following statements hold.
\begin{enumerate}[(a)]
\item
Let $t\geq2$ and $n\geq2t^{2}+2t$ be two integers.
If $m>\binom{n-1}{2}+t-1$, then $G$ is $t$-tough.
\item
Let $n\geq\max\{(\frac{2t^{2}+3t+1}{t})\delta+\frac{1}{2t}+5,\frac{1}{6}(\frac{\delta(\delta+4)+1}{t}+9\delta+6)\}$, where $\frac{1}{t}\geq1$ and $\delta\geq t+1$ are integers.
If $m>\binom{n-\frac{\delta}{t}-1}{2}+\frac{(\delta+t)\delta}{t}$, then $G$ is $t$-tough.
\end{enumerate}
\end{lemma}

Now, we divide the Theorem \ref{th:1.3.} into two parts to proof it.\\

\noindent\textbf{Proof of  Theorem \ref{th:1.3.}.(a).}

Since $Tr(v)\geq d(v)+2(n-1-d(v))=2(n-1)-d(v)$ with equality holds if and only if the maximum distance between $v$ and other vertices in $G$ is at most 2.

So
$$W(G)=\frac{1}{2}\sum\limits_{v\in V(G)}Tr(v)\geq\frac{1}{2}\sum\limits_{v\in V(G)}[2(n-1)-d(v)]=n(n-1)-m$$
with equality holds if and only if the maximum distance between $v$ and other vertices in $G$ is at most 2.

By Lemma \ref{le:2.3.}, we get $\eta_{1}(G)\geq\frac{4W(G)}{n}\geq4(n-1)-\frac{4m}{n}$.
Combining with the condition of Theorem \ref{th:1.3.}, we obtain
$$4(n-1)-\frac{4m}{n}\leq\eta_{1}(G)<\frac{2n^{2}+2n-4t}{n},$$
which implies $m>\binom{n-1}{2}+t-1$.
According to Lemma \ref{le:4.1.}, $G$ is $t$-tough.

This completes the proof.
$\hfill\square$\\

\noindent\textbf{Proof of Theorem \ref{th:1.3.}.(b).}

By the above analyses, we know that $\eta_{1}(G)\geq\frac{4W(G)}{n}\geq4(n-1)-\frac{4m}{n}$.
Combining with the condition of Theorem \ref{th:1.3.}, we obtain
$$4(n-1)-\frac{4m}{n}\leq\eta_{1}(G)<2n+\frac{4\delta}{t}+2-\frac{2(\delta+t)(2\delta t+\delta+2t)}{nt^{2}},$$
which implies $m>\binom{n-\frac{\delta}{t}-1}{2}+\frac{(\delta+t)\delta}{t}$.
According to Lemma \ref{le:4.1.}, $G$ is $t$-tough.

This completes the proof.
$\hfill\square$\\

\end{document}